\documentclass[aap,MSNbibl,seceqn,dvips]{arximspdf}

\doi{10.1214/13-AAP930} 
\volume{24}
\issue{2}
\pubyear{2014}
\firstpage{599}
\lastpage{615}

\makeatletter
\newcommand{\eqref}[1]{(\ref{#1})}
\newcommand{\last}{\mathrm{last}}
\newtheorem{theorem}{Theorem}[section]
\newtheorem{lemma}[theorem]{Lemma}
\newtheorem{proposition}[theorem]{Proposition}
\makeatother

\begin{document}
\begin{frontmatter}

\title{Path properties of the disordered pinning model in~the delocalized regime}
\runtitle{Path delocalization}

\begin{aug}
\author[A]{\fnms{Kenneth S.} \snm{Alexander}\thanksref{t1,t3}\ead[label=e1]{alexandr@usc.edu}}
\and
\author[B]{\fnms{Nikos} \snm{Zygouras}\corref{}\thanksref{t2,t3}\ead[label=e2]{N.Zygouras@warwick.ac.uk}}
\runauthor{K. S. Alexander and N. Zygouras}
\affiliation{University of Southern California and University of
Warwick}
\address[A]{Department of Mathematics\\
University of Southern California\\
Los Angeles, California 90089-2532\\
USA\\
\printead{e1}}
\address[B]{Department of Statistics\\
University of Warwick\\
Coventry CV4 7AL\\
United Kingdom\\
\printead{e2}}
\end{aug}
\thankstext{t1}{Supported by NSF Grant DMS-08-04934.}
\thankstext{t2}{Supported by a Marie Curie International Reintegration Grant within the 7th European Community Framework Programme, IRG-246809.}
\thankstext{t3}{Supported in part by the North America Travel Fund of the University of Warwick.}

\received{\smonth{10} \syear{2012}}
\revised{\smonth{3} \syear{2013}}

%
\begin{abstract}
We study the path properties of a random polymer attracted to a defect
line by a potential with disorder, and we prove that in the delocalized
regime, at any temperature, the number of contacts with the defect line
remains in a certain sense ``tight in probability'' as the polymer
length varies. On the other hand we show that at sufficiently low
temperature, there exists a.s. a subsequence where the number of
contacts grows like the log of the length of the polymer.
\end{abstract}

%

\end{frontmatter}

\section{Introduction}
The disordered pinning model has attracted significant attention in
recent years. One reason is that it is one of the very few models where
the effect of disorder on the critical properties can be identified
with large precision. In particular, there exists a fairly satisfactory
knowledge on whether and how much the critical point, which separates
its localized and delocalized regime, changes under the presence of
disorder \cite{A,GLT2}. Furthermore, the mechanism that defines it is
present in multiple physical models, and therefore it provides a step
to understand the effect of disorder in more complicated systems---we
refer to the recent monograph \cite{G} for related references.

Before going into detail let us define the model. We first consider a
sequence of i.i.d. variables $(\omega_n)_{n\in\mathbb{Z}}$, which play
the role of disorder. The assumptions on this sequence are in general
mild, for example, mean zero and exponential moments. We denote the
joint distribution of this sequence by $\mathbb{P}$. The model involves
also a renewal sequence $(\tau_n)_{n\in\mathbb{N}}$ on
$\mathbb{N}=\{0,1,2,\ldots\}$, that is, a point process such that the
gaps (or interarrival times) $\sigma_n:=\tau_{n+1}-\tau_n$ are
independent and identically distributed. This renewal process should be
viewed physically as the set of contact points with $\{0\}\times
\mathbb{N}$ of the space-time trajectory of a Markov process
$(X_n)_{n\in\mathbb{N}}$ whose state space contains a designated site
0, with this trajectory representing the spatial configuration of the
polymer. Since the interaction between the Markov process and the
disorder comes only at contact times with $\{0\} \times\mathbb{N} $,
the only relevant information is the renewal sequence
$\tau=(\tau_n)_{n\in\mathbb{N}}$, consisting of the contact points of
the path $(X_n)_{n\in\mathbb{N}}$ with $\{0\} \times\mathbb{N} $.
Therefore we only need to define the statistics of this renewal
process, whose law we will denote by $P$. In particular, we define
$\tau_0=0$ and assume that for some $\alpha\geq0$ and slowly varying
function~$\phi(n)$,
\[
K(n):=P(\tau_1=n)=\frac{\phi(n)}{n^{1+\alpha}},\qquad n\geq1.
\]
We will assume that $\sum_{n\geq1}K(n)=1$, that is, that the renewal
is recurrent. We will~also need the quantity $K^+(l) = \sum_{n>l} K(n)$.

The polymer measure can now be defined by
\[
dP_{n,\omega}^{\beta,h}:=\frac{1}{Z^{\beta,h}_{n,\omega
}}e^{\mathcal{H}_{n,\omega}^{\beta,u}}\,dP,
\]
where $\mathcal{H}_{n,\omega}^{\beta,u}:=\sum_{i=0}^n(\beta\omega
_i+h)\delta_i$ and $\delta_i=1_{i\in\tau}$.
The partition function $ Z^{\beta,h}_{n,\omega }$ is defined by
\[
Z^{\beta,h}_{n,\omega } = E \bigl[e^{\mathcal{H}_{N,\omega}^{\beta,h}}
\bigr].
\]
The polymer measure rewards paths for which the $\omega_i$ values are
large at the times of renewals. It will also be useful to consider the
constrained polymer measure
\[
dP_{n,\omega}^{\beta,h,c}:=\frac{1}{Z^{\beta,h,c}_{n,\omega
}}e^{\mathcal{H}_{n,\omega}^{\beta,u}}
\delta_n\,dP,
\]
where we restrict the polymer to have a renewal at time $n$. Here the
constrained partition function is
\[
Z^{\beta,h,c}_{n,\omega } = E \bigl[e^{\mathcal{H}_{N,\omega
}^{\beta,h}}\delta_n
\bigr].
\]
More generally for a collection $A$ of trajectories we define
\[
Z^{\beta,h}_{n,\omega }(A) = E \bigl[e^{\mathcal{H}_{N,\omega
}^{\beta,h}}; A \bigr].
\]
We will also need the notation
\[
Z^{\beta,h}_{[m,n],\omega }=Z^{\beta,h}_{n-m,\theta_{m}\omega },
\]
where $n\geq m$ and $\theta_m\omega (i)=\omega (i+m)$, for
$i=1,2,\ldots\,$.

As already mentioned, the pinning polymer exhibits a nontrivial
localization/delocalization transition, which is often
quantified via the strict positivity of the free energy. To be more
precise, let us define the quenched free energy of the pinning
polymer to be the $\mathbb{P}$-a.s. limit
\[
f_q(\beta,h)=\lim_{n\to\infty}\frac{1}{n}\log
Z^{\beta,h}_{n,\omega }.
\]
We refer the reader to \cite{G}, Chapter~3, for the existence of this
limit. The localized regime is defined as
\[
\mathcal{L}= \bigl\{(\beta,h)\dvtx f_q(\beta,h)>0 \bigr\},
\]
and the delocalized regime as
\[
\mathcal{D}= \bigl\{(\beta,h)\dvtx f_q(\beta,h)=0 \bigr\}.
\]
The free energy is monotone in $h$ so the two regimes are separated by
a critical line and we can define the quenched critical point
$h_c(\beta)$ as
\[
h_c(\beta)=\sup \bigl\{h\dvtx f_q(\beta,h)=0 \bigr\}.
\]
Let $M(\beta)=\mathbb{E}[e^{\beta\omega_1}]$ be the moment
generating function of $\omega_1$. For the corresponding annealed
model, with partition function $\mathbb{E}Z^{\beta,h}_{n,\omega }$
and free energy
\[
f_a(\beta,h)=\lim_{n\to\infty}\frac{1}{n}\log
\mathbb{E}Z^{\beta,h}_{n,\omega },
\]
the corresponding critical point is
%
\begin{equation}
\label{hcan} h_c^{\mathrm{ann}}(\beta) = -\log M(\beta).
\end{equation}
The question of the path behavior of the quenched model for
$h<h_c(\beta)$ is of particular interest when $h_c^{\mathrm{ann}}(\beta
)<h_c(\beta)$, so we summarize what has been proved about such an inequality.
It is known from~\cite{A,T} (for Gaussian disorder) and
from~\cite{L} (for general disorder) that for small $\beta$, $h_c(\beta
)=h_c^{\mathrm{ann}}(\beta)$, for $\alpha<1/2$ as well as for $\alpha=1/2$
and $\sum_{n\geq1}(n\phi(n)^2)^{-1}<\infty$. On the other hand,
from \mbox{\cite{A,AZ09,DGLT,AS},} for Gaussian
disorder, for $1/2<\alpha< 1$, there exists a constant $c$ and a
slowly varying function $\psi$ related to $\phi$ and $\alpha$ such
that for all small $\beta$,
\[
c^{-1}\beta^{2\alpha/(2\alpha-1)}\psi \biggl(\frac{1}{\beta
}
\biggr)<h_c(\beta)-h_c^{\mathrm{ann}}(\beta)< c
\beta^{2\alpha/(2\alpha-1)}\psi \biggl(\frac{1}{\beta} \biggr),
\]
while for $\alpha=1$,
\[
c^{-1}\beta^2\psi \biggl(\frac{1}{\beta} \biggr) <
h_c(\beta )-h_c^{\mathrm{ann}}(\beta). %
\]
A matching upper bound is also expected to hold but has not been proved.
For $\alpha> 1$,
\[
c^{-1}\beta^2<h_c(\beta)-h_c^{\mathrm{ann}}(
\beta)< c\beta^2. %
\]
The case $\alpha=1/2$ is marginal and not fully understood. It is
believed that $h_c(\beta)>h_c^{\mathrm{ann}}(\beta)$ for every $\beta$, as
long as $\sum_n 1/(n\phi(n)^2)=\infty$. This inequality has been
confirmed under some stronger hypotheses in \cite{AZ09,GLT1},
for Gaussian disorder, and (most nearly optimally, for general
disorder) in \cite{GLT2}.
For all $\alpha>0$, for large $\beta$ the critical points are shown
in \cite{To} to be distinct provided the disorder is unbounded, but
for $\alpha=0$ they are equal for all $\beta>0$ \cite{AZ10}. Theorem
1.5 of \cite{CdH} shows that for $\alpha>1/2$ the critical points are
different for all values of $\beta>0$.

The use of the terms localization/delocalization can be understood
better by relating the quenched free energy to the portion of time the
polymer spends on the defect line $\{0\}\times\mathbb{N}$. In
particular, from \cite{GT}, $f_q(\beta,\cdot)$ is differentiable for
all $h\neq h_c(\beta)$ with
\[
\frac{d}{dh}f_q(\beta,h)=\lim_{n\to\infty}
E^{P_{n,\omega}^{\beta,h}} \Biggl[\frac{1}{n}\sum_{i=1}^n
\delta_i \Biggr],
\]
and therefore we can interpret the localized regime as the regime where
the polymer spends a positive fraction of time on the defect line,
while in the delocalized regime it spends a zero fraction of time on
the defect line. While this is quite satisfactory in the localized
regime, and further detailed studies on the path properties in the
localized regime have been made in \cite{GT}, it provides a rather
incomplete picture in the delocalized one---it only allows one to
conclude that the number of contacts is $o(n)$. It was proven in
\cite{GT2} that the number of contacts is at most of order $\log n$ in
the delocalized regime. This was actually done for the related
copolymer model, but its extension to the pinning model is
straightforward \cite{G}. More precisely, for every $h<h_c(\beta)$,
there exists a constant $C_{\beta,h}$ such that
\[
\limsup_{n\to\infty} \mathbb{E}P_{n,\omega}^{\beta,h}
\bigl(\bigl|\tau \cap[1,n]\bigr|>C_{\beta,h}\log n \bigr)=0.
\]
This result was further extended to an a.s. statement in \cite{M}: for
$h<h_c(\beta)$ and for every $C>(1+\alpha)/(h_c(\beta)-h)$, we have
\[
\limsup_{n\to\infty} P_{n,\omega}^{\beta,h} \bigl(\bigl|\tau\cap
[1,n]\bigr|>C\log n \bigr)=0,\qquad\mathbb{P}\mbox{-a.s.}
\]
By analogy to the homogeneous pinning model (see \cite{G}, Chapter~8),
one might expect that the number of contacts with the defect line
should remain bounded in the whole delocalized regime. Nevertheless,
the picture has been unclear in the disordered case, since stretches of
unusual disorder values could typically attract the polymer back to the
defect line a number of times growing to infinity with $n$. The open
questions are discussed in \cite{G}, Section~8.5.
In this work we clarify and complete the picture for behavior in
probability. In fact, we will prove a stronger result, namely, that the
last contact of the polymer happens at distance $O(1)$ from the origin.
In particular, let
\[
\tau_{\last}=\max\{j\leq n\dvtx\delta_j=1 \}.
\]
We then have the following theorem.

\begin{theorem} \label{delocalization}
Suppose $\alpha>0$, $\sum_n K(n)=1$ and that $\omega_1$ has exponential
moments of all orders. For all $\beta$, $\varepsilon >0$ and for all
$h<h_c(\beta)$ we have that
\[
\limsup_{N\to\infty}\limsup_{n\to\infty}\mathbb{P}
\bigl(P_{n,\omega}^{\beta,h} (\tau_{\last}>N )>\varepsilon
\bigr)=0.
\]
\end{theorem}

One may ask whether this can be made an almost-sure result for
$h<h_c(\beta)$, of the form
\[
\limsup_{N\to\infty}\limsup_{n\to\infty}P_{n,\omega}^{\beta,h}
(\tau_{\last}>N )=0, \qquad\mathbb{P}\mbox{-a.s.},
\]
or if the number of contacts is a.s. finite, that is,
\[
\limsup_{N\to\infty}\limsup_{n\to\infty}P_{n,\omega}^{\beta,h}
\bigl(\bigl|\tau\cap[0,n]\bigr|>N \bigr)=0, \qquad\mathbb{P}\mbox{-a.s.}
\]
The next theorem shows that the answer is no, at least for large $\beta
$. Instead, for $h$ between $h_c^{\mathrm{ann}}(\beta)$ and $h_c(\beta)$,
infinitely often as $n\to\infty$, there will be an exceptionally rich
segment of $\omega $ near $n$, which will (with high $P_{n,\omega
}^{\beta,h}$-probability) induce the polymer to come to 0 and then
make a number of returns of order $\log n$. For $t>0$ let
%
\begin{equation}
\label{htdef} h_t(\beta):= -(1+t\alpha)\log M \biggl(
\frac{\beta}{1+t\alpha} \biggr).
\end{equation}
Since $\log M$ is nondecreasing and convex on $[0,\infty)$ with $\log
M(0)=0$, it is easy to see that $h_t(\beta)$ is nondecreasing in $t$
for fixed $\beta$. Recall \eqref{hcan}; by \cite{Ca},
equation~(3.7), for all $\beta>0$ we have
%
\begin{equation}
\label{hcbeta} -\log M(\beta)=h_c^{\mathrm{ann}}(
\beta)=h_0(\beta) \leq h_c(\beta) \leq h_1(
\beta).
\end{equation}
By \cite{To}, Theorem 3.1, given $0<\varepsilon<1$, for large $\beta$
we have
%
\begin{equation}
\label{hcbeta2} h_c(\beta) > h_{1-\varepsilon}(\beta).
\end{equation}
We are now ready to state our second main result.

\begin{theorem} \label{longreturns}
Suppose $\omega $ is unbounded with all exponential moments finite.
Given $\varepsilon >0$, there exists $\beta_0(\varepsilon)$ and $\nu
(\beta,h)>0$ such that for
\[
\beta>\beta_0 \quad\mbox{and} \quad h>h_\varepsilon(\beta),
\]
we have
\[
\limsup_{n\to\infty}P_{n,\omega}^{\beta,h} \bigl(\bigl|\tau\cap
[0,n]\bigr|>\nu\log n \bigr)=1, \qquad\mathbb{P}\mbox{-a.s.}
\]
\end{theorem}

By \eqref{hcbeta2}, Theorem \ref{longreturns} with $\varepsilon <1/2$
includes at least the interval of values $h\in[h_\varepsilon(\beta
),h_c(\beta)]$ below $h_c(\beta)$, which in turn (for large $\beta$)
includes $h\in[h_\varepsilon(\beta),h_{1-\varepsilon}(\beta)]$. The
path behavior in the regime of Theorem \ref{longreturns} is therefore
in contrast with that for $h<h_c^{\mathrm{ann}}(\beta)$, where, in fact, the
number of contacts remains tight for the measures averaged over the
disorder; see \cite{GT2}, Remark 1.5.

The next two sections are devoted to the proofs of each theorem, respectively.

\section{\texorpdfstring{Proof of Theorem \protect\ref{delocalization}}
{Proof of Theorem 1.1}}

It will be convenient to introduce generic constants. Specifically, $C$
will denote a generic constant whose value might be different in
different appearances. If we want to distinguish between constants we
will enumerate them, for example,~$C_1$,~$C_2$, etc. When we want to
emphasize the dependence of a generic constant on some parameters, we
will include the symbols of these parameters as a subscript. In
particular, we use the notation $C_\alpha$ for a generic constant which
will depend on the parameter $\alpha$ and the slowly varying function
$\phi$ of the renewal process. To simplify the notation we will also
defer from using the integer part $[x]$ and simply write $x$, which
should not lead to any confusion in the contexts where we use it. Let
us define the events
\[
E_{n,N}= \bigl\{\bigl|\tau\cap[0,n]\bigr|>N \bigr\}, \qquad E_{[m,n],N}=
\bigl\{\bigl|\tau\cap [m,n]\bigr|>N \bigr\}.
\]
%

In proving Theorem \ref{delocalization} we will make use of the
following theorem, which was proved in \cite{M}.

\begin{theorem}[(\cite{M})] \label{Mourrat}
Let $\beta\geq0$ and $h<h_c(\beta)$. Then:
\begin{longlist}[(iii)]
\item[(i)] For $\mathbb{P}$-a.e. environment $\omega $, we have
\[
\sum_{n=0}^\infty Z^{\beta,h,c}_{n,\omega }<+
\infty.
\]

\item[(ii)] For every $\varepsilon >0$ and for $\mathbb{P}$-a.e.
environment $\omega $, there exists $N_\varepsilon (\omega)>0$
such that for all $N\geq N_\varepsilon $, we have that
\[
\sum_{n=0}^\infty Z^{\beta,h,c}_{n,\omega }(E_{n,N})
\leq\sum_{k=N}^\infty e^{-k(h_c(\beta)-h-\varepsilon )}.
\]

\item[(iii)] For every constant $C>\frac{1+\alpha}{h_c(\beta)-h}$
and for $\mathbb{P}$-a.e. environment $\omega $, we have
\[
P_{n,\omega}^{\beta,h,c}(E_{n,C\log n})\longrightarrow0,\qquad\mbox{as
} n\to\infty.
\]
\end{longlist}
\end{theorem}

The quantity $\mathcal{Z}(\omega )=\sum_{n=0}^\infty Z^{\beta,h,c}_{n,\omega }$, which is a.s. finite, will play an important role,
as will the reversed process $\mathcal{Z}_n(\omega)=\sum_{m=-\infty}^n
Z^{\beta,h,c}_{[m,n],\omega }$, which for any fixed $n$ has the same
distribution as $\mathcal{Z}(\omega )$. Note that we think here of the
polymer path starting at point $n$ and going backwards in time, which
is why we have defined the disorder on the whole of~$\mathbb{Z}$.

Here is a sketch of the proof. The event $\{\tau_{\last}>N\}$ is
contained in the union of the following events, where $C_1$, $b>0$ are
constants, with $b$ small:
\begin{longlist}[(a)]
\item[(a)] there are more than $C_1\log n$ returns by time $n$;

\item[(b)] there are fewer than $C_1\log n$ returns, and no gap
between returns exceeds~$bn$;
\item[(c)] $\tau_{\last}>N$, there are fewer than $C_1\log n$
returns, and some gap between returns inside $[0,n]$ exceeds
$bn$;

\item[(d)] $\tau_{\last}>N$, there are fewer than $C_1\log n$
returns, and the incomplete gap $[\tau_{\last},n]$ exceeds
$bn$.
\end{longlist}
The Gibbs probability of (a) can be controlled by a variant of Theorem
\ref{Mourrat}(iii), (b)~can be controlled using the small probability
of the event under the free measure and (d) is relatively
straightforward, so the main work is (c). The segment to the left of
the size-$bn$ gap corresponds to a term in the sum $\mathcal {Z}(\omega
)$, and (after we ``tie down'' the right end of the polymer by adding a
visit at time $n$) the segment to the right corresponds to a term in
$\mathcal{Z}_n(\omega )$, so we make use of Theorem~\ref{Mourrat}(i)
and a bound for the probability of a big gap under the free measure.

Let us make note here of the trivial lower bound
%
\begin{equation}
\label{onejump} Z^{\beta,h}_{n,\omega } \geq K^+(n)e^{\beta\omega_0+h},
\end{equation}
which comes from the trajectory having no renewals after time 0.

We will need the following analog of Theorem \ref{Mourrat}(iii), for
the free polymer measure.

\begin{lemma}\label{Muncon}
Let $\beta\geq0$ and $h<h_c(\beta)$. Then for all $C_1>\frac{\alpha
}{h_c(\beta)-h}$ and for $\mathbb{P}$-a.e. environment $\omega $, we
have
\[
P_{n,\omega}^{\beta,h}(E_{n,C_1\log n})\longrightarrow0,\qquad\mbox{as
} n\to\infty.
\]
\end{lemma}

\begin{pf}
Let $\varepsilon >0$ satisfy $C_1>\frac{\alpha+\varepsilon }{h_c(\beta
)-h-\varepsilon }$. Using Theorem \ref{Mourrat}(ii) and
\eqref{onejump}, for some $C_2=C_2(\beta,h,\varepsilon,\alpha)$, we
have for large $n$
%
\begin{eqnarray}
\label{numbounda} Z^{\beta,h}_{n,\omega }(E_{n,C_1\log n}) &=& \sum
_{j=1}^n Z^{\beta,h,c}_{j,\omega }(E_{j,C_1\log n})
K^+(n-j)
\nonumber
\\
&\leq&\sum_{k=C_1\log n}^\infty e^{-k(h_c(\beta)-h-\varepsilon )}
\nonumber
\\
&\leq& C_2n^{-(\alpha+\varepsilon )}
\\
&\leq& C_2 K^+(n)e^{\beta\omega_0+h}n^{-\varepsilon /2}
\nonumber
\\
&\leq& C_2 n^{-\varepsilon /2} Z^{\beta,h}_{n,\omega },
\nonumber
\end{eqnarray}
and the lemma follows.
\end{pf}

Proposition \ref{nojump} below will show that the probability is small
for having fewer than $C_1\log n$ renewals without some gap $\sigma$
exceeding $bn$, when $b$ is chosen sufficiently small. Let us denote
this gap event by $A_{b,n}$; more precisely, let
%
\begin{eqnarray}\label{adef}
A_{b,n}'&=& \bigl\{\mbox{$\tau$: there exist $i,j\in[0,n]$, $j-i\geq bn$,}\nonumber
\\
&&\hspace*{49pt}\mbox{ such that } \tau\cap[i,j]=\{i,j\} \bigr\}
\\
A_{b,n}''&=& \{ \tau\dvtx
\tau_{\last} \leq n-bn \}\nonumber
\\
A_{b,n} &=& A_{b,n}' \cup A_{b,n}''.
\nonumber
\end{eqnarray}
We first prove an analogous statement for the free renewal process.

\begin{lemma}\label{decaybig}
Given $C_1$ as in Lemma \ref{Muncon} and $b \in(0,1/2)$, for
sufficiently large~$n$, we have
\[
P \bigl(E_{n,C_1\log n}^{c}\cap A_{b,n}^c
\bigr)\leq n^{-\alpha/9b}.
\]
\end{lemma}

\begin{pf}
When the event $E_{n,C_1\log n}^{c}\cap A_{b,n}^c$ occurs, there exists
$l \leq C_1 \log n$ such that
\[
\sigma_1+\cdots+\sigma_l = \tau_{\last}
\in(n-bn,n] \quad\mbox{and}\quad\max_{i \leq l} \sigma_i
< bn.
\]
Among these first $l$ jumps, the total length of all jumps having
individual length $\sigma_i \leq n/4C_1\log n$ is at most $n/4$, so the
total length of all jumps with individual length $\sigma_i \in
[n/4C_1\log n,bn)$ is at least $n/4$. This means there must be at least
$1/4b$ values $\sigma_i \geq n/4C_1\log n$ among $\sigma_1,\ldots,\sigma_{C_1\log n}$. Presuming $n$ is large, we have
\[
p_n:= P \biggl(\sigma_1 \geq\frac{n}{4C_1\log n} \biggr)
\leq n^{-\alpha/2}.
\]
Let $k_n$ be the integer part of $C_1\log n$, and let $r$ be the least
integer greater than or equal to $1/4b$. Then for large $n$,
\begin{eqnarray*}
P \bigl(E_{n,C_1\log n}^{c}\cap A_{b,n}^c
\bigr) &\leq& P \biggl( \biggl\vert \biggl\{ i \leq k_n\dvtx
\sigma_i \geq\frac
{n}{4C_1\log n} \biggr\} \biggr\vert \geq r \biggr)
\\
&\leq& \pmatrix{k_n \cr r} p_n^r
\leq(k_np_n)^r \leq n^{-\alpha/9b}.
\end{eqnarray*}\upqed
\end{pf}

\begin{proposition} \label{nojump}
Given $C_1>0$ as in Lemma \ref{Muncon} and given $\beta,h$, for $b>0$
sufficiently small,
\[
P_{n,\omega}^{\beta,h} \bigl(E_{n,C_1\log n}^{c}\cap
A_{b,n}^c \bigr) \to0 \qquad\mbox{a.s. as } n \to\infty.
\]
\end{proposition}

\begin{pf}
We have from Lemma \ref{decaybig} that if $b$ is sufficiently small,
then for large~$n$,
%
\begin{eqnarray}
\label{meandecay} && \frac{1}{K^+(n)} \mathbb{E} \bigl[ Z^{\beta,h}
\bigl(E_{n,C_1\log
n}^{c}\cap A_{b,n}^c \bigr)
\bigr]
\nonumber
\\
&&\qquad \leq \frac{C_\alpha n^\alpha}{\phi(n)} e^{(\log M(\beta)+h)C_1\log
n} P \bigl(E_{n,C_1\log n}^{c}
\cap A_{b,n}^c \bigr)
\\
&&\qquad \leq \frac{1}{n^3}.
\nonumber
\end{eqnarray}
Therefore for all $\eta>0$,
%
\begin{eqnarray}
\label{io1} && \mathbb{P} \bigl(P_{n,\omega}^{\beta,h}
\bigl(E_{n,C_1\log
n}^{c}\cap A_{b,n}^c \bigr) >
\eta\mbox{ i.o.} \bigr)
\nonumber
\\[3pt]
&&\qquad \leq\mathbb{P} \bigl(Z^{\beta,h}_{n,\omega }
\bigl(E_{n,C_1\log
n}^{c}\cap A_{b,n}^c \bigr) >
\eta K^+(n)e^{\beta\omega_0+h} \mbox{ i.o.} \bigr)
\nonumber
\\[-5pt]
\\[-5pt]
&& \qquad \leq\mathbb{P} \bigl(Z^{\beta,h}_{n,\omega }
\bigl(E_{n,C_1\log
n}^{c}\cap A_{b,n}^c \bigr) >
\eta K^+(n)n^{-1} \mbox{ i.o.} \bigr)
\nonumber
\\[3pt]
&&\quad\qquad{} + \mathbb{P} \biggl(e^{\beta\omega_0+h} < \frac{1}{n} \mbox{
i.o.} \biggr).
\nonumber
\end{eqnarray}
Now the second probability on the right-hand side of \eqref{io1} is 0,
and by \eqref{meandecay}, for the first probability on the right-hand
side, we have
%
\begin{eqnarray}
&& \mathbb{P} \bigl(Z^{\beta,h}_{n,\omega } \bigl(E_{n,C_1\log
n}^{c}
\cap A_{b,n}^c \bigr) > \eta K^+(n)n^{-1} \bigr)
\nonumber
\\[3pt]
&&\qquad \leq\frac{n}{\eta K^+(n)} \mathbb{E} \bigl[ Z^{\beta,h}
\bigl(E_{n,C_1\log n}^{c}\cap A_{b,n}^c \bigr)
\bigr]
\\[3pt]
&&\qquad \leq\frac{1}{\eta n^2}.
\nonumber
\end{eqnarray}
Summing over $n$ and applying the Borel--Cantelli lemma completes the
proof.
\end{pf}

The next proposition, together with Lemma \ref{Muncon} and Proposition
\ref{nojump}, shows that with probability tending to one, the first
big gap, of length at least $bn$, brings the polymer out of $[0,n]$.

\begin{proposition} \label{bigjump}
For every $b,\varepsilon >0$ we have
%
\begin{equation}
\label{limsups} \lim_{n\to\infty}\mathbb{P} \bigl(P_{n,\omega}^{\beta,h}
\bigl(A_{b,n}' \bigr) > \varepsilon \bigr)=0.
\end{equation}
\end{proposition}

\begin{pf}
Let $0<\theta<1$. Then, summing over possible locations $[n_1,n_2]$
for the interval of the first long jump, we have
%
\begin{eqnarray}\label{prop25eq1}
Z^{\beta,h}_{n,\omega} \bigl(A_{b,n}'
\bigr) &\leq& \sum_{n_1}\sum
_{n_1+bn<n_2\leq n} Z^{\beta,h,c}_{n_1,\omega} K(n_2-n_1)
Z^{\beta,h}_{[n_2,n],\omega}
\nonumber
\\[3pt]
&=& \sum_{n_1=0}^n\sum
_{\max(n_1+bn,n-n^\theta)<n_2\leq n} Z^{\beta,h,c}_{n_1,\omega}K(n_2-n_1)
Z^{\beta,h}_{[n_2,n],\omega}
\\[3pt]
&&{} + \sum_{n_1=0}^n\sum
_{n_1+bn <n_2\leq n-n^\theta} Z^{\beta,h,c}_{n_1,\omega}K(n_2-n_1)
Z^{\beta,h}_{[n_2,n],\omega}.
\nonumber
\end{eqnarray}
Using \eqref{onejump}, we can bound the first term on the right-hand
side of \eqref{prop25eq1} by
%
\begin{eqnarray}
\label{prop25eq2} \hspace*{20pt}&& \sum_{n_1}\sum
_{\max(n_1+bn,n-n^\theta)<n_2\leq n}\ \sum_{l\leq n-n_2}
Z^{\beta,h,c}_{n_1,\omega}\,K(n_2-n_1)\,Z^{\beta,h,c}_{[n_2,n-l],\omega}\,K^+(l)
\nonumber
\\
&&\qquad \leq C K(bn) \sum_{n_1}\sum
_{\max(n_1+bn,n-n^\theta)<n_2\leq n}\ \sum_{l\leq n-n_2}
Z^{\beta,h,c}_{n_1,\omega}\,Z^{\beta,h,c}_{[n_2,n-l],\omega}\,K^+(l)
\nonumber
\\
&&\qquad \leq \frac{C}{b^{1+\alpha}} K(n) n^\theta \sum
_{n_1}\sum_{\max(n_1+bn,n-n^\theta)<n_2\leq n}\ \sum
_{l\leq n-n_2} Z^{\beta,h,c}_{n_1,\omega}\,Z^{\beta,h,c}_{[n_2,n-l],\omega}\,K(l)
\\
&&\qquad \leq \frac{C}{b^{1+\alpha}} K^+(n) n^{\theta-1} e^{-(\beta\omega_n+h)}\sum
_{n_1}\sum_{\max(n_1+bn,n-n^\theta)<n_2\leq n}
Z^{\beta,h,c}_{n_1,\omega}\,Z^{\beta,h,c}_{[n_2,n],\omega}
\nonumber
\\
&&\qquad \leq \frac{C}{b^{1+\alpha}} Z^{\beta,h}_{n,\omega} n^{\theta-1}
e^{-(\beta\omega_0+h)}e^{-(\beta\omega_n+h)} \mathcal{Z}(\omega) \mathcal{Z}_n(
\omega).\nonumber
\end{eqnarray}
The second term on the right-hand side of \eqref{prop25eq1} is bounded
by
%
\begin{eqnarray}
\label{prop25eq3} \qquad&& \sum_{n_1}\sum
_{n_1+bn <n_2<n-n^\theta}\ \sum_{l\leq n-n_2}
Z^{\beta,h,c}_{n_1,\omega }\,K(n_2-n_1)
\,Z^{\beta,h,c}_{[n_2,n-l],\omega }\,K^+(l)
\nonumber
\\
&&\qquad \leq C K(bn)\sum_{n_1}\sum
_{n_1+bn <n_2\leq n-n^\theta}\ \sum_{l\leq n-n_2}
Z^{\beta,h,c}_{n_1,\omega }\,Z^{\beta,h,c}_{[n_2,n-l],\omega }\,K^+(l)
\nonumber
\\
&&\qquad \leq \frac{C}{b^{1+\alpha}} K(n) n \sum_{n_1}
\sum_{n_1+bn
<n_2\leq n-n^\theta}\ \sum_{l\leq n-n_2}
Z^{\beta,h,c}_{n_1,\omega
}\,Z^{\beta,h,c}_{[n_2,n-l],\omega }\,K(l)
\nonumber
\\[-8pt]
\\[-8pt]
&&\qquad \leq \frac{C}{b^{1+\alpha}} K^+(n) e^{-(\beta\omega_n+h)} \sum
_{n_1}\sum_{n_1+bn <n_2\leq n-n^\theta}
Z^{\beta,h,c}_{n_1,\omega }\,Z^{\beta,h,c}_{[n_2,n],\omega }
\nonumber
\\
&&\qquad \leq \frac{C}{b^{1+\alpha}} Z^{\beta,h}_{n,\omega}
e^{-(\beta\omega_0+h)}e^{-(\beta\omega_n+h)} \sum_{n_1=0}^\infty
Z^{\beta,h,c}_{n_1,\omega } \sum_{n_2=-\infty}^{n-n^\theta}
Z^{\beta,h,c}_{[n_2,n],\omega }
\nonumber
\\
&&\qquad \leq \frac{C}{b^{1+\alpha}} Z^{\beta,h}_{n,\omega}
e^{-(\beta\omega_0+h)}e^{-(\beta\omega_n+h)} \mathcal{Z}(\omega ) \sum
_{n_2=-\infty}^{n-n^\theta} Z^{\beta,h,c}_{[n_2,n],\omega }.
\nonumber
\end{eqnarray}
From \eqref{prop25eq1}, \eqref{prop25eq2} and \eqref{prop25eq3} we
have that
%
\begin{eqnarray}
\label{twoterms} P_{n,\omega}^{\beta,h} \bigl(A_{b,n}'
\bigr)&\leq& \frac{C}{b^{1+\alpha}} n^{\theta-1} e^{-(\beta\omega
_0+h)}e^{-(\beta\omega_n+h)}
\mathcal{Z}(\omega) \mathcal {Z}_n(\omega)
\nonumber
\\[-9pt]
\\[-9pt]
&&{} + \frac{C}{b^{1+\alpha}} e^{-(\beta\omega_0+h)}e^{-(\beta
\omega_n+h)} \mathcal{Z}(\omega )
\sum_{n_2=-\infty
}^{n-n^\theta} Z^{\beta,h,c}_{[n_2,n],\omega }.
\nonumber
\end{eqnarray}
Now $\mathcal{Z}(\omega )$ and $\mathcal{Z}_n(\omega )$ are finite
almost surely and equidistributed, so the first term on the right in
\eqref{twoterms} converges to 0 in $\mathbb{P}$-probability. The sum
on the right-hand side of \eqref{twoterms} has the same distribution as
\[
\sum_{m=n^\theta}^\infty Z^{\beta,h,c}_{m,\omega },
\]
so by Theorem \ref{Mourrat}(i), it converges to 0 in probability.
Hence the second term on the right-hand side of \eqref{twoterms} also
converges to 0 in probability, and the proof is complete.
\end{pf}

We can now complete the proof of our first theorem.
\begin{pf*}{Proof of Theorem \ref{delocalization}}
For $b>0$ we have
\begin{eqnarray*}
P_{n,\omega}^{\beta,h}(\tau_{\last}>N) &\leq&
P_{n,\omega}^{\beta,h} \bigl(E_{n,C_1\log n}^c\cap
A_{b,n}^c \bigr)+ P_{n,\omega}^{\beta,h}(E_{n,C_1\log n})
\\[-2pt]
&&{}+
P_{n,\omega}^{\beta,h} \bigl(A_{b,n}' \bigr) +P_{n,\omega}^{\beta,h} \bigl(\{\tau_{\last}>N \cap
A_{b,n}'' \bigr).
\end{eqnarray*}
By Proposition \ref{nojump}, with the choice of sufficiently small
$b>0$, and Lemma \ref{Muncon}, respectively, we have that the first and
second terms in the above expression converge to zero $\mathbb
{P}$-a.s., while by Proposition \ref{bigjump}, the third term converges
to 0 in $\mathbb{P}$-probability. Therefore, it only remains to check
that
\[
\limsup_{N\to\infty}\limsup_{n\to\infty}\mathbb{P}
\bigl(P_{n,\omega}^{\beta,h} \bigl(\{\tau_{\last}>N\} \cap
A_{b,n}'' \bigr) >\varepsilon \bigr)=0.
\]
To this end we have
\begin{eqnarray*}
&& \mathbb{P} \bigl(P_{n,\omega}^{\beta,h} \bigl(\{\tau_{\last}
\geq N\} \cap A_{b,n}'' \bigr) >\varepsilon
\bigr)
\\[-2pt]
&&\qquad \leq \mathbb{P} \bigl(Z^{\beta,h}_{n,\omega } \bigl(\{
\tau_{\last}\geq N\} \cap A_{b,n}''
\bigr) >\varepsilon K^+(n) e^{\beta\omega_0+h} \bigr)
\\[-2pt]
&&\qquad =\mathbb{P} \biggl(\sum_{N\leq n_1<n-bn}Z^{\beta,h,c}_{n_1,\omega
}K^+(n-n_1)
>\varepsilon K^+(n) e^{\beta\omega_0+h} \biggr)
\\[-2pt]
&&\qquad \leq \mathbb{P} \biggl(\sum_{N\leq n_1<n-bn}Z^{\beta,h,c}_{n_1,\omega }
>\varepsilon C_{\alpha,b}\,e^{\beta\omega_0+h} \biggr)
\\[-2pt]
&&\qquad \leq\mathbb{P} \biggl(\sum_{n_1\geq N}Z^{\beta,h,c}_{n_1,\omega }
>\varepsilon C_{\alpha,b}\,e^{\beta\omega_0+h} \biggr),
\end{eqnarray*}
and by Theorem \ref{Mourrat}(i), the latter tends to 0 as $N\to\infty$.
\end{pf*}

The analog of Theorem \ref{delocalization} also holds for the
constrained case, that is, for~$P_{n,\omega}^{\beta,h,c}$, in the sense
that the rightmost contact point in $[0,\frac{n}{2}] $ and the leftmost
contact point in $[\frac{n}{2},n]$ occur at distances $O(1)$ from 0 and
$n$, respectively. To quantify things, let us denote
\[
\hat\tau_{\last}=\max \biggl\{j \in \biggl[0, \frac{n}{2} \biggr]
\dvtx\delta_j=1 \biggr\}
\]
and
\[
\check\tau_{\last}=\min \biggl\{j \in \biggl[\frac{n}{2},n \biggr]
\dvtx\delta_j=1 \biggr\}.
\]
Then we have the following.

\begin{theorem}\label{consdel}
Suppose $\alpha>0$, $\sum_n K(n)=1$ and that $\omega_1$ has exponential
moments of all orders. For all $\beta,\varepsilon,\delta
>0$ and for all $h<h_c(\beta)$ there exist $n_0(\varepsilon,\delta
)$, $N_0(\varepsilon,\delta)$ and $M_0(\varepsilon,\delta)$, such
that for all $n>n_0(\varepsilon,\delta)$, $N>N_0(\varepsilon,\delta)$, $M>M_0(\varepsilon,\delta)$
\[
\mathbb{P} \bigl(P_{n,\omega}^{\beta,h,c} \bigl(\{\hat\tau
_{\last}>N\}\cup\{\check\tau_{\last}<n-M\} \bigr)>\varepsilon
\bigr)<\delta.
\]
\end{theorem}

\begin{pf}
Notice that in the constrained case $A_{b,n}=A_{b,n}'$, and we have
\begin{eqnarray*}
&& P_{n,\omega}^{\beta,h,c} \bigl(\{\hat\tau_{\last}>N\}\cup\{
\check \tau_{\last}<n-M\} \bigr)
\\
&&\qquad \leq P_{n,\omega}^{\beta,h,c} \bigl(E_{n,C_1\log n}^c
\cap A_{b,n}^c \bigr)+ P^{\beta,h,c}_{n,\omega} (E_{n,C_1\log n})
\\
&&\qquad\quad{} + P_{n,\omega}^{\beta,h,c} \bigl( \bigl(\{\hat
\tau_{\last}>N\} \cup\{\check\tau_{\last}<n-M\} \bigr) \cap
A_{b,n}' \bigr).
\end{eqnarray*}
By a straightforward modification of Proposition \ref{nojump},
Theorem~\ref{Mourrat}(iii) and Lemma~\ref{Muncon}, the first two terms converge
to zero as $n$ tends to infinity, once $b$ is chosen small enough.
Regarding the third term, notice that by symmetry it is sufficient to
control $Z^{\beta,h,c}_{n,\omega } (\{\hat\tau_{\last}>N\} \cap
A_{b,n}' )$. We make two sums according to whether the first big gap
$[n_1,n_2]$ ends before or after the midpoint $n/2$. Specifically, we
have
\begin{eqnarray*}
&& Z^{\beta,h,c}_{n,\omega
} \bigl(\{\hat\tau_{\last}>N\} \cap
A_{b,n}' \bigr)
\\
&&\qquad \leq \sum_{n_1>N} \sum_{\max(n_1+bn,n/2)<n_2\leq n} Z^{\beta,h,c}_{n_1,\omega}
K(n_2-n_1) Z^{\beta,h,c}_{[n_2,n],\omega }
\\
&&\quad\qquad{}+ \sum_{n_1}\sum_{n_1+bn<n_2\leq n/2}
Z^{\beta,h,c}_{n_1,\omega} K(n_2-n_1)
Z^{\beta,h,c}_{[n_2,n],\omega }.
\end{eqnarray*}
Following the same (and actually more direct) steps as in the proof of
Proposition~\ref{bigjump} we can bound the above by
\[
C_b e^{-(\beta\omega_0+h)}e^{-(\beta\omega_n+h)}Z^{\beta,h,c}_{n,\omega
} \biggl(\sum_{n_1>N} Z^{\beta,h}_{n_1,\omega } \mathcal{Z}_n(\omega )
+\mathcal{Z}(\omega )\sum _{n_2<n/2}Z^{\beta ,h,c}_{[n_2,n],\omega }
\biggr),
\]
and the rest follows as in Proposition \ref{bigjump}.
\end{pf}

\section{\texorpdfstring{Proof of Theorem \protect\ref{longreturns}}
{Proof of Theorem 1.2}}
We begin again with a sketch. Assume for simplicity that $K(1)>0$.
Suppose there is a ``rich segment'' of $[0,n]$ of length at least
$\gamma\log n$ in which the average of the disorder is at least $u$;
here $\gamma$ is small and $u$ is large. (We show that such a rich
segment exists for infinitely many $n$.) We consider the contribution
to the partition function from two different sets of trajectories:
\begin{longlist}[(a)]
\item[(a)] the single trajectory which returns at every site of the
rich segment, and nowhere else;

\item[(b)] those trajectories which make at most $\nu\log n$ returns,
with $\nu$ small.
\end{longlist}
We show that (up to slowly varying correction factors) the contribution
from (a) is at least a certain inverse power $n^{-\alpha+\kappa}$,
while a.s., except for finitely many $n$, the contribution from (b) is
bounded by the smaller inverse power $n^{-\alpha+\kappa/2}$. The
Gibbs probability of (b) is bounded by the ratio of the two
contributions, hence by~$n^{-\kappa/2}$, so it approaches~0.

We will need the following lemma, which is an elementary fact about
convex functions.

\begin{lemma} \label{convex}
Suppose $\Psi$ is nondecreasing and convex on $[0,\infty)$ with $\Psi
(0)=0$ and $\Psi'(x) \to\infty$ as $x \to\infty$. Then for all $s>1$,
\[
\Psi(sx) - s\Psi(x) \to\infty\qquad\mbox{as } x \to\infty.
\]
\end{lemma}

\begin{pf}
Since $\Psi'$ is nondecreasing, for $s>1$ and $x>1$, we have
%
\begin{eqnarray}
\label{derives} \Psi(sx) - s\Psi(x) &=& (s-1) \int_0^x
\bigl( \Psi' \bigl(x+(s-1)t \bigr) - \Psi'(t) \bigr)
\, dt
\\
&\geq& (s-1) \int_0^1 \bigl(
\Psi' \bigl(x+(s-1)t \bigr) - \Psi'(t) \bigr)\, dt
\\
&\geq& (s-1) \bigl(\Psi'(x) - \Psi'(1) \bigr)
\\
&\to& \infty \qquad \mbox{as } x \to \infty.
\end{eqnarray}
Here the first inequality follows from the fact that the integrand in
nonnegative.\vadjust{\goodbreak}
\end{pf}

\begin{pf*}{Proof of Theorem \ref{longreturns}}
Recall the definition of $h_t(\beta)$ from \eqref{htdef}. Suppose
$h=h_t(\beta)$ with $t>\varepsilon $. If $t\geq1$, then $h \geq
h_c(\beta)$ by \eqref{hcbeta}, so we need only consider $t<1$. Let $r =
\min\{j\dvtx K(j)>0\}$, let $\gamma,u>0$ to be specified, define
\[
J_n = \{n-ir\dvtx 0\leq i \leq\gamma\log n - 1\}, \qquad \overline{
\omega}_{J_n} = \frac{1}{|J_n|}\sum_{j\in J_n}
\omega_j,
\]
and define the event
\[
D_n^{u,\gamma} = \{\omega \dvtx \overline{\omega}_{J_n}
\geq u\}.
\]
We can bound $Z^{\beta,h}_{n,\omega }$ below by the contribution from
the path which makes returns precisely at the times in $J_n$, obtaining
that for large $n$, for all \mbox{$\omega \in D_n^{u,\gamma}$},
%
\begin{eqnarray}
\label{Zlower} Z^{\beta,h}_{n,\omega } &\geq& e^{(\beta u + h)|J_n|} K(n-
\gamma r\log n)K(r)^{|J_n|-1}
\nonumber
\\[-8pt]
\\[-8pt]
&\geq& \frac{1}{2}n^{-(1+\alpha)}\phi(n) \exp \biggl(\gamma \biggl(\beta
u + h - \log\frac{1}{K(r)} \biggr)\log n \biggr).\nonumber
\end{eqnarray}
Let $\Phi$ be the large deviation rate function related to $\omega $,
and let $\delta>0$ to be specified. For large $n$ we have
%
\begin{equation}
\label{largedev} \mathbb{P} \bigl(D_n^{u,\gamma} \bigr) \geq
e^{-(1+\delta)\Phi(u)\gamma\log n}.
\end{equation}
Since all exponential moments of $\omega$ are finite, we have $\Phi
(u)/u\to\infty$ as $u\to\infty$.
Recalling that $\log M(\beta) = \sup\{\beta u - \Phi(u)\dvtx u\in
\mathbb{R}\}$, we can therefore choose \mbox{$u=u_\beta$} to satisfy
\[
\beta u - \Phi(u) = \log M(\beta).
\]
For $\beta$ sufficiently large (depending on $\varepsilon $), since
$\Phi'(u) \to\infty$ as $u\to\infty$, we have by Lemma \ref{convex}
that
%
\begin{equation}
\label{Krcost} \log M(\beta) - (1+\varepsilon \alpha)\log M \biggl(
\frac{\beta
}{1+\varepsilon \alpha} \biggr) > \log\frac{1}{K(r)}
\end{equation}
or equivalently,
\[
\beta u_\beta+ h_{\varepsilon }(\beta) - \log\frac{1}{K(r)} >
\Phi (u_\beta).
\]

We now choose $\delta$ to satisfy
\[
\beta u_\beta+ h_{\varepsilon }(\beta) - \log\frac{1}{K(r)} > (1+
\delta)\Phi(u_\beta)
\]
and then $\gamma$ to satisfy
%
\begin{equation}
\label{gamma} \beta u_\beta+ h_{\varepsilon }(\beta) - \log
\frac{1}{K(r)} > \frac
{1}{\gamma} > (1+\delta)\Phi(u_\beta).
\end{equation}
Define $\kappa>0$ by
%
\begin{eqnarray}
\label{kappadef} \gamma \biggl(\beta u_\beta+ h_{\varepsilon }(\beta)
- \log\frac
{1}{K(r)} \biggr) &=& 1+\kappa,
\end{eqnarray}
so that by \eqref{Zlower}, for all $\omega \in D_n^{u,\gamma}$,
%
\begin{equation}
\label{kappabound} Z^{\beta,h}_{n,\omega } \geq\tfrac{1}{2}n^{-\alpha+\kappa}
\phi(n).
\end{equation}
We select a subsequence of the events $\{D_n^{u,\gamma}\}$ that are
independent, as follows. Fix $n_0$ and given $n_0,\ldots,n_j$ define
$n_{j+1} = n_j + 2r\gamma\log n_j$. Then $n_j \sim2r\gamma j\log j$ as
$j\to\infty$, and it is easily checked that, provided $n_0$ is
sufficiently large, the events $\{D_{n_j}^{u,\gamma}, j\geq0\}$ are
independent. With \eqref{largedev} and \eqref{gamma} this shows that
%
\begin{equation}
\label{io} \sum_j \mathbb{P}
\bigl(D_{n_j}^{u,\gamma} \bigr) = \infty\quad\mbox{so}\quad
\mathbb{P} \bigl(D_n^{u,\gamma} \mbox{ i.o.} \bigr) = 1.
\end{equation}
Let us now choose
%
\begin{equation}
\label{mlambdanu} m > \frac{4}{\kappa}, \qquad\lambda=2 \biggl(
\frac{1}{m}\log M(m\beta)+h \biggr), \qquad\nu= \frac{\kappa}{2\lambda},
\end{equation}
with $m$ an integer. We claim that
%
\begin{equation}
\label{polybound} \mathbb{P} \bigl(Z^{\beta,h}_{n,\omega }
\bigl(E_{n,\nu\log n}^c \bigr) > n^{-\alpha+\lambda\nu}\phi(n) \mbox{
i.o.} \bigr) =0.
\end{equation}
This is plausible because for appropriate $\lambda$, $\nu\log n$
visits should not likely yield more than $\lambda\nu\log n$ energy
above the ``immediate escape'' value, which is approximately the log of
$K^+(n)$, that is, approximately $-\alpha\log n$.
Assuming this claim, we use \eqref{io} to conclude that
\[
\mathbb{P} \bigl(D_n \cap \bigl\{ Z^{\beta,h}_{n,\omega }
\bigl(E_{n,\nu\log n}^c \bigr) < n^{-\alpha+\lambda\nu}\phi(n) \bigr\}
\mbox{ i.o.} \bigr) = 1,
\]
which with \eqref{kappabound} shows that
\[
\mathbb{P} \bigl(P_{n,\omega}^{\beta,h} \bigl(E_{n,\nu\log n}^c
\bigr) < 2n^{-\kappa/2} \mbox{ i.o.} \bigr)=1,
\]
which proves the theorem.

It remains to prove \eqref{polybound}. Observe that by Chebyshev's
inequality we have
%
\begin{eqnarray}\label{polybound2}
&& \mathbb{P} \bigl(Z^{\beta,h}_{n,\omega }
\bigl(E_{n,\nu\log
n}^c \bigr) > n^{-\alpha+\lambda\nu}\phi(n) \bigr)
\nonumber\\[-8pt]\\[-8pt]
&&\qquad \leq
\bigl(n^{-\alpha}\phi(n) \bigr)^{-m} n^{-m\lambda\nu} \mathbb{E}
\bigl[ \bigl(Z_{n,\omega }^{\beta,h} \bigl(E_{n,\nu\log n}^c
\bigr) \bigr)^m \bigr].\nonumber
\end{eqnarray}
Denoting by $E^{\otimes m}$ the expectation over $m$ independent copies
of the renewal $\tau$, we see
that the expectation on the right-hand side of \eqref{polybound2} can
be written as
\[
E^{\otimes m} \bigl[ e^{\sum_{i=1}^n (\log M(\beta(\delta
_i^{(1)}+\cdots+\delta_i^{(m)})) +h(\delta_i^{(1)}+\cdots+\delta
_i^{(m)}) )}; \bigl(E_{n,\nu\log n}^c
\bigr)^{\otimes m} \bigr],
\]
where $(E_{n,\nu\log n}^c)^{\otimes m}$ is the $m$-fold product of
$E_{n,\nu\log n}^c$. Using the convexity of $\log M(\beta)$ we have
\[
\log M(\beta k) \leq\frac{k}{m} \log M(\beta m) \qquad\mbox{for all } k
\leq m,
\]
so we can bound the above expectation by
%
\begin{eqnarray}
\label{mfoldmean} && E^{\otimes m} \bigl[ e^{\sum_{i=1}^n ((1/m)\log M(m\beta
) +h )(\delta_i^{(1)}+\cdots+\delta_i^{(m)})};
\bigl(E_{n,\nu\log
n}^c \bigr)^{\otimes m} \bigr]
\nonumber
\\[-8pt]
\\[-8pt]
&& \qquad< e^{ ((1/m)\log M(m\beta) +h ) m\nu\log n} P \bigl(E_{n,\nu
\log n}^c
\bigr)^m.
\nonumber
\end{eqnarray}
We use $A_{b,n}$ from \eqref{adef}. By Lemma \ref{decaybig} we have
for $b$ sufficiently small and then $n$ sufficiently large,
\begin{eqnarray}
P \bigl(E_{n,\nu\log n}^c \bigr) &\leq& P \bigl(E_{n,\nu\log n}^c
\cap A_{b,n}^c \bigr) + \sum_{j=1}^{\nu\log n}
P(\sigma_j >bn)
\nonumber
\\
&\leq& n^{-2\alpha} + \nu K^+(bn)\log n
\nonumber
\\
&\leq& C_b\nu\log n \frac{\phi(n)}{n^\alpha}.
\nonumber
\end{eqnarray}
Inserting this into \eqref{mfoldmean} and the result into \eqref
{polybound2}, and considering our choice of $\lambda,m,\nu$, we
obtain that
\[
\mathbb{P} \bigl(Z^{\beta,h}_{n,\omega } \bigl(E_{n,\nu\log n}^c
\bigr) > n^{-\alpha+\lambda\nu}\phi(n) \bigr) \leq (C_b\nu\log
n)^m n^{-m\kappa/4},
\]
which, by the choice of $m$ in \eqref{mlambdanu} and the
Borel--Cantelli lemma, completes the proof.
\end{pf*}

If we do not assume $\beta$ large in Theorem \ref{longreturns}, then in
the proof, the entropy cost $\log1/K(r)$ per visit to $J_n$ will not be
exceeded by the energy gain; in more concrete terms, \eqref {Krcost}
will fail. The entropy cost can be reduced by visiting only a small
fraction of the sites in an interval of form $[n-\gamma\log n,n]$, but
then the interval length $\gamma\log n$ (where the disorder average
exceeds $u_\beta$) must be much larger than in the large-$\beta$ proof,
reducing the probability of such an interval. It is not clear whether
there is a strategy (in place of the present ``visit all sites of
$J_n$'') of sufficiently low entropy cost so that the interval of large
average disorder values can be exploited, and therefore it seems
unclear whether a variant of Theorem \ref{longreturns} should be true
for small~$\beta$.\looseness=1


\section*{Acknowledgements}

The authors would like to thank the referee for improved proofs of
Lemmas \ref{decaybig} and \ref{convex}.
Part of this work was
done during the School and Conference on Random Polymers and Related
Topics, at the National University of Singapore. The authors would like
to thank the institute for the hospitality.



\printaddresses

\end{document}